\documentclass[leqno,12pt]{amsart}
\usepackage{amsfonts}
\usepackage{amsmath,amssymb,amsthm}

\newcommand{\R}{\mathbb R}

\newcommand{\E}{\mathbb E}

\renewcommand{\span}{\mathrm{span}}
\newcommand{\tr}{\mathrm{tr}}

\newtheorem{thm}{Theorem}[section]

\theoremstyle{definition}

\theoremstyle{remark}

\newcommand{\ds}{\displaystyle}

\begin{document}

\title[Meridian Surfaces  with 1-type Gauss map in Minkowski 4-Space] {Meridian Surfaces of Elliptic or Hyperbolic Type with Pointwise 1-type Gauss Map in Minkowski 4-Space}

\author{Kadri Arslan and Velichka Milousheva}
\address{Uluda\u{g} University, Art and Science Faculty, Department of Mathematics, 16059 Bursa,
Turkey} \email{arslan@uludag.edu.tr}
\address{Institute of Mathematics and Informatics, Bulgarian Academy of Sciences,
Acad. G. Bonchev Str. bl. 8, 1113, Sofia, Bulgaria;   "L. Karavelov"
Civil Engineering Higher School, 175 Suhodolska Str., 1373 Sofia,
Bulgaria} \email{vmil@math.bas.bg}


\begin{abstract}
In the present paper we consider a special class of spacelike
surfaces in the Minkowski 4-space which are one-parameter systems
of meridians of the rotational hypersurface with timelike or
spacelike axis. They are called meridian surfaces of elliptic or
hyperbolic type, respectively. We study these surfaces with
respect to their Gauss map. We find all meridian surfaces of
elliptic or hyperbolic type with harmonic Gauss map and give the
complete classification of  meridian surfaces of elliptic or
hyperbolic type with  pointwise 1-type Gauss map.
\end{abstract}

\maketitle

\section{Introduction}

The study of submanifolds of Euclidean space or pseudo-Euclidean
space via the notion of finite type immersions began in the late
1970's with the papers \cite{Ch1,Ch2} of B.-Y. Chen and has been
extensively carried out since then. An isometric immersion $x:M$
$\rightarrow $ $\E^{m}$ of a submanifold $M$ in Euclidean
$m$-space $\E^{m}$  (or  pseudo-Euclidean space $\E^m_s$)  is said
to be of \emph{finite type} \cite{Ch1}, if $x$ identified with the
position vector field of $M$ in $\E^{m}$ (or $\E^m_s$) can be
expressed as a finite sum of eigenvectors of the Laplacian $\Delta
$ of $M$, i.e.
\begin{equation*}
x=x_{0}+\sum_{i=1}^{k}x_{i},
\end{equation*}
where $x_{0}$ is a constant map, $x_{1},x_{2},...,x_{k}$ are
non-constant maps such that $\Delta x_i=\lambda _{i}x_{i},$
$\lambda _{i}\in \mathbb{R}$, $1\leq i\leq k.$ If $\lambda
_{1},\lambda _{2},...,\lambda _{k}$ are different, then $M$ is
said to be of \emph{$k$-type}. Many results on finite type
immersions have been collected in the survey paper \cite{Ch3}.

The notion of finite type immersion is naturally extended to the
Gauss map $G$ on $M$ by  B.-Y. Chen and  P. Piccinni \cite{CP}.
Thus,  a submanifold $M$ of an  Euclidean (or pseudo-Euclidean
space) is said to have  \emph{1-type Gauss map} $G$, if $G$
satisfies $\Delta G=a (G+C)$ for some $a \in \mathbb{R}$ and some
constant vector $C$ (see, for example, \cite{BB}, \cite{BCV},
\cite{BV}, \cite{KY1}).

However, the Laplacian of the Gauss map of some well-known
surfaces such as the helicoid, the catenoid, and the right cone in
the Euclidean 3-space $\mathbb{E}^{3}$, the helicoids of 1st, 2nd,
and 3rd kind, conjugate of Enneper's surface of 2nd kind and
B-scrolls in the Minkowski 3-space $\E^3_1$,  the generalized
catenoids, and Enneper's hypersurfaces in $\E^{n+1}_1$  takes a
somewhat different form, namely, $\Delta G=\lambda (G+C)$ for some
non-constant smooth function $\lambda$ and some constant vector
$C$.  Therefore, it is worth studying the class of surfaces
satisfying such an equation.

We use the following definition:  a submanifold $M$ of the
Euclidean space $\E^m$  (or pseudo-Euclidean space $\E^m_s$) is
said to have \emph{pointwise 1-type Gauss map} if its Gauss map
$G$ satisfies
\begin{equation} \notag
\Delta G=\lambda (G+C)
\end{equation}
for some non-zero smooth function $\lambda $ on $M$ and some
constant vector $C$. A pointwise 1-type Gauss map is called
\emph{proper} if the function $\lambda $  is non-constant. A
submanifold with pointwise 1-type Gauss map is said to be of \emph{first kind} if the vector $C$ is zero. Otherwise, it is said
to be of \emph{second kind} \cite{CCK}.

Classification results on surfaces with pointwise 1-type Gauss map
in Minkowski space have been obtained in the last few years. For
example, in \cite{KY3} Y.  Kim and  D. Yoon studied ruled surfaces with
1-type Gauss map in Minkowski space $\E^m_1$ and gave a complete
classification of null scrolls with 1-type Gauss map. The
classification of ruled surfaces with pointwise 1-type Gauss map
of first kind in Minkowski space $\E^3_1$ is given in
\cite{KY2}. Ruled surfaces with pointwise 1-type Gauss map  of
second kind in Minkowski 3-space were classified in \cite{CKY}.
The complete classification of flat rotation surfaces with  pointwise
1-type Gauss map  in the 4-dimensional pseudo-Euclidean space
$\E^4_2$ is given in \cite{KY2-a}.

\vskip 2mm Basic source of examples of surfaces in the
four-dimensional Euclidean or Minkowski space are the
meridian surfaces. Meridian surfaces in the Euclidean 4-space
$\R^4$ are defined in \cite{GM2} as special class of surfaces,
which are one-parameter systems of meridians of the standard
rotational hypersurface in $\R^4$. In \cite{ABM}  we studied the
meridian surfaces with pointwise 1-type Gauss map. We showed that
a meridian surface in $\R^4$  has a harmonic Gauss map if and only
if it is part of a plane. We gave necessary and sufficient
conditions for a meridian surface to have pointwise 1-type Gauss
map and found all meridian surfaces with pointwise 1-type Gauss
map of first and second kind. The meridian surfaces of Weingarten type
are described in \cite{ABBO}.

The meridian surfaces of elliptic or hyperbolic type in the
Minkowski 4-space $\R^4_1$ are constructed in \cite{GM6} similarly
to the Euclidean case.  They are
 two-dimensional  spacelike surfaces in $\R^4_1$   which are
one-parameter systems of meridians of the rotational hypersurface
with timelike or spacelike axis, respectively. Recently, some
special classes of meridian surfaces of elliptic or hyperbolic
type have been  classified. For example, marginally trapped
meridian surfaces of elliptic or hyperbolic type are described in
\cite{GM6}. The complete classification of meridian surfaces of
elliptic or hyperbolic type with constant Gauss curvature or with
constant mean curvature is given in \cite{GM-new}. The
 Chen meridian surfaces and  the meridian surfaces with parallel
normal bundle are also classified in \cite{GM-new}.

In the present paper we study meridian surfaces of elliptic or
hyperbolic type in $\R^4_1$ with respect to their Gauss map. In
Theorem  \ref{T:harmonic-ell} and Theorem \ref{T:harmonic-hyp} we describe all meridian surfaces of elliptic or
hyperbolic type with harmonic Gauss map. We give the complete
classification of meridian surfaces of elliptic or hyperbolic type
with  pointwise 1-type Gauss map of first  kind in
Theorem \ref{T:First kind - ell} and Theorem \ref{T:First kind - hyp}, respectively.
The meridian surfaces of elliptic or hyperbolic type
with  pointwise 1-type Gauss map of second  kind are classified in Theorem \ref{T:Second kind - ell} and Theorem \ref{T:Second kind - hyp}, respectively.

\section{Preliminaries} \label{S:Pre}

Let $\R^4_1$  be the four-dimensional Minkowski space  endowed
with the metric $\langle , \rangle$ of signature $(3,1)$ and
$Oe_1e_2e_3e_4$ be a fixed orthonormal coordinate system such
that $\langle e_1, e_1 \rangle  = \langle e_2, e_2 \rangle  =
\langle e_3, e_3 \rangle  = 1, \, \langle e_4, e_4 \rangle  = -1$. The standard flat metric is
given in local coordinates by $dx_1^2 + dx_2^2 + dx_3^2 -dx_4^2.$

A surface $M$
in $\R^4_1$ is said to be
\emph{spacelike} if $\langle , \rangle$ induces  a Riemannian
metric $g$ on $M$. Thus at each point $p$ of a spacelike surface
$M$ we have the following decomposition:
$$\R^4_1 = T_pM \oplus N_pM$$
with the property that the restriction of the metric $\langle ,
\rangle$ onto the tangent space $T_pM$ is of signature $(2,0)$,
and the restriction of the metric $\langle , \rangle$ onto the
normal space $N_pM$ is of signature $(1,1)$.

We denote by $\nabla'$ and $\nabla$ the Levi Civita connections on $\R^4_1$ and $M$, respectively.
Let $x$ and $y$ be vector fields tangent to $M$ and $\xi$ be a normal vector field.
The formulas of Gauss and Weingarten giving the decompositions of the vector fields $\nabla'_xy$ and
$\nabla'_x \xi$ into tangent and normal components are given, respectively, by
$$\begin{array}{l}
\vspace{2mm}
\nabla'_xy = \nabla_xy + \sigma(x,y);\\
\vspace{2mm}
\nabla'_x \xi = - A_{\xi} x + D_x \xi,
\end{array}$$
where  $\sigma$ is the second fundamental tensor,  $D$ is the normal connection,
and $A_{\xi}$ is the shape operator  with respect to $\xi$.

The mean curvature vector  field $H$ of $M$ is defined as
$H = \ds{\frac{1}{2}\,  \tr\, \sigma}$.
A submanifold $M$ is said to be minimal (respectively, totally geodesic) if $%
H = 0$ (respectively, $\sigma = 0$).
A surface $M$ in the Minkowski 4-space is called \emph{marginally trapped} \cite{Chen2}, if its mean curvature vector field $H$ is lightlike at each point, i.e.
$H \neq 0$, $\langle H, H \rangle =0$.

\vskip 1mm The Gauss map $G$ of a submanifold $M$ of
$\E^m$ is defined as follows. Let $G(n,m)$ be the Grassmannian
manifold consisting of all oriented $n$-planes through the origin
of $\mathbb{E}^{m}$ and $\wedge ^{n}\mathbb{E}^{m}$ be the vector
space obtained by the exterior product of $n$ vectors in
$\mathbb{E}^{m}$. In a natural way, we can identify $\wedge
^{n}\mathbb{E}^{m}$ with the Euclidean space $\mathbb{E}^{N}$,
where $N=\left(
\begin{array}{c}
m \\
n
\end{array}
\right)$.
 Let $\left\{ e_{1},...,e_{n},e_{n+1},\dots,e_{m}\right\} $ be a
 local orthonormal frame field in $\mathbb{E}^{m}$ such that $e_{1},e_{2},\dots,$ $e_{n}$ are tangent to $M$ and
 $e_{n+1},e_{n+2},\dots,e_{m}$ are
 normal to $M$.
The map $G:M\rightarrow G(n,m)$ defined by $%
G(p)=(e_{1}\wedge e_{2}\wedge \dots \wedge $ $e_{n})(p)$ is called the \emph{Gauss
map} of $M$. It is a smooth map which carries a point $p$ in $M$ into the
oriented $n$-plane in $\mathbb{E}^{m}$ obtained by the parallel translation
of the tangent space of $M$ at $p$ in $\mathbb{E}^{m}$.

In a similar way one can consider the Gauss map of  a submanifold $M$ of pseudo-Euclidean space  $\E^m_s$.

For any  function $f$ on $M$ the Laplacian of $f$ is given by the formula

\begin{equation}\notag
\Delta f =-\sum_{i}(\nabla'_{e_{i}}\nabla'
_{e_{i}}f -\nabla'_{\nabla _{e_{i}}e_{i}}f ),
\end{equation}
where $\nabla'$ is the Levi-Civita connection of $\E^m$ (or $\E^m_s$)  and $\nabla$ is the induced connection on $M$.

\section{Meridian surfaces of elliptic or hyperbolic type} \label{S:Mer}

Meridian surfaces in the Minkowski space $\E^4_1$ are special families of two-dimensional  spacelike surfaces lying
 on rotational hypersurfaces in $\R^4_1$  with
timelike or spacelike axis, which are constructed as follows.

Let $f = f(u), \,\, g = g(u)$ be smooth functions, defined in an
interval $I \subset \R$, such that $(f'(u))^2 - (g'(u))^2 > 0,
\,\, u \in I$. We assume that $f(u)>0, \,\, u \in I$. The standard
rotational hypersurface $\mathcal{M}'$ in $\R^4_1$, obtained by
the rotation of the meridian curve $m: u \rightarrow (f(u), g(u))$
about the $Oe_4$-axis,  is parameterized as follows:
$$\mathcal{M}': Z(u,w^1,w^2) = f(u)\, \cos w^1 \cos w^2 \,e_1 +  f(u)\, \cos w^1 \sin w^2 \,e_2 + f(u)\, \sin w^1 \,e_3 + g(u) \,e_4.$$
The rotational hypersurface $\mathcal{M}'$ is a two-parameter system of meridians. If $w^1 =
w^1(v)$, $w^2=w^2(v), \,\, v \in J, \, J \subset \R$,  we can consider the two-dimensional
 surface $\mathcal{M}'_m$ lying on $\mathcal{M}'$, constructed  in the following way:
\begin{equation}  \notag
\mathcal{M}'_m: z(u,v) = Z(u,w^1(v),w^2(v)), \quad u \in I, \, v \in J.
\end{equation}
Since $\mathcal{M}'_m$ is a
one-parameter system of meridians of $\mathcal{M}'$, it is called
a  \emph{meridian surface of elliptic type}  \cite{GM6}.

If we denote $l(w^1,w^2) = \cos w^1 \cos w^2 \,e_1 + \cos w^1 \sin w^2 \,e_2 + \sin w^1 \,e_3$, then the surface $\mathcal{M}'_m$ is parameterized by
\begin{equation} \label{E:Eq2}
\mathcal{M}'_m: z(u,v) = f(u) \, l(v) + g(u)\, e_4, \quad u \in I, \, v \in J.
\end{equation}
Note that $l(w^1,w^2)$ is  the unit position vector of the 2-dimensional sphere $S^2(1)$ lying in the Euclidean space $\R^3 = \span \{e_1, e_2, e_3\}$
and centered at the origin $O$.

We assume that the smooth curve  $c: l = l(v) = l(w^1(v),w^2(v)), \, v \in J$  on
$S^2(1)$ is parameterized by the arc-length, i.e. $\langle l'(v), l'(v) \rangle = 1$. Let
 $t(v) = l'(v)$  be the tangent vector field of $c$. Since $\langle t(v), t(v) \rangle = 1$,
 $\langle l(v), l(v) \rangle = 1$, and $\langle t(v), l(v) \rangle = 0$,
 there exists a unique (up to a sign)
 vector field $n(v)$, such that
$\{ l(v), t(v), n(v)\}$ is an orthonormal frame field in $\R^3$. With respect to this
 frame field we have the following Frenet formulas of $c$ on $S^2(1)$:
\begin{equation} \label{E:Eq3}
\begin{array}{l}
\vspace{2mm}
l' = t;\\
\vspace{2mm}
t' = \kappa \,n - l;\\
\vspace{2mm} n' = - \kappa \,t,
\end{array}
\end{equation}
 where $\kappa (v)= \langle t'(v), n(v) \rangle$ is the spherical curvature of $c$.

Without loss of generality we assume that $(f'(u))^2 - (g'(u))^2 = 1$.
The tangent space of $\mathcal{M}'_m$ is spanned by the vector fields:
$$z_u = f' \,l + g'\,e_4; \qquad  z_v = f\,t,$$
so, the coefficients of the first fundamental form of $\mathcal{M}'_m$  are $E = 1; \, F = 0; \, G = f^2(u) >0$.
Hence, the first fundamental form is positive
definite, i.e. $\mathcal{M}'_m$ is a spacelike surface.

Denote  $x = z_u,\,\, y = \ds{\frac{z_v}{f} = t}$ and
consider the following orthonormal normal frame field:
$$n_1 = n(v); \qquad n_2 = g'(u)\,l(v) + f'(u) \, e_4.$$
Thus we obtain a frame field $\{x,y, n_1, n_2\}$ of $\mathcal{M}'_m$, such that $\langle n_1, n_1 \rangle =1$,
$\langle n_2, n_2 \rangle =- 1$, $\langle n_1, n_2 \rangle =0$.

Taking into account \eqref{E:Eq3} we get
the following derivative formulas \cite{GM-new}:
\begin{equation} \label{E:Eq5}
\begin{array}{ll}
\vspace{2mm} \nabla'_xx = \kappa_m\,n_2; & \qquad
\nabla'_x n_1 = 0;\\
\vspace{2mm} \nabla'_xy = 0;  & \qquad
\nabla'_y n_1 = \ds{ - \frac{\kappa}{f}\,y};\\
\vspace{2mm} \nabla'_yx = \ds{\frac{f'}{f}}\,y;  & \qquad
\nabla'_x n_2 = \kappa_m \,x;\\
\vspace{2mm} \nabla'_yy = \ds{- \frac{f'}{f}\,x + \frac{\kappa}{f}\,n_1 + \frac{g'}{f} \, n_2}; & \qquad
\nabla'_y n_2 = \ds{ \frac{g'}{f}\,y},
\end{array}
\end{equation}
where $\kappa_m (u)= f'(u) g''(u) - g'(u) f''(u)$ is the curvature of the
meridian curve $m$, and  $\kappa = \kappa (v)$ is the spherical curvature of $c$.

\vskip 2mm
In a similar way one can consider meridian surfaces lying on the
 rotational hypersurface in $\R^4_1$ with spacelike axis.
Let $f = f(u), \,\, g = g(u)$ be smooth functions, defined in an
interval $I \subset \R$, such that $(f'(u))^2 + (g'(u))^2 >0$, $f(u)>0, \,\, u \in I$.
The rotational hypersurface $\mathcal{M}''$ in $\R^4_1$, obtained
by the rotation of the meridian curve $m: u \rightarrow (f(u),
g(u))$ about the $Oe_1$-axis  is parameterized as follows:
$$\mathcal{M}'': Z(u,w^1,w^2) = g(u) \,e_1 +  f(u)\, \cosh w^1 \cos w^2 \,e_2 +  f(u)\, \cosh w^1 \sin w^2 \,e_3+ f(u)\, \sinh w^1 \,e_4.$$
If $w^1 = w^1(v), \, w^2=w^2(v), \,\, v \in J, \,J \subset
\R$, we consider the surface $\mathcal{M}''_m$ in $\R^4_1$ defined by
\begin{equation}  \notag
\mathcal{M}''_m: z(u,v) = Z(u,w^1(v),w^2(v)),\quad u \in I, \, v \in J.
\end{equation}
 $\mathcal{M}''_m$ is a one-parameter system of
meridians of $\mathcal{M}''$ and is called a
\emph{meridian surface of hyperbolic type} \cite{GM6}.

If we denote $l(w^1,w^2) = \cosh w^1 \cos w^2 \,e_2 +  \cosh w^1 \sin w^2 \,e_3 + \sinh w^1 \,e_4$, then the surface $\mathcal{M}''_m$
is given by
\begin{equation} \label{E:Eq6}
\mathcal{M}''_m: z(u,v) = f(u) \, l(v) + g(u)\, e_1, \quad u \in I, \, v \in J,
\end{equation}
$l(w^1,w^2)$ being the unit position vector of the  timelike sphere $S^2_1(1)$  in the Minkowski space $\R^3_1 = \span \{e_2, e_3, e_4\}$,
i.e. $S^2_1(1) = \{ V \in \R^3_1: \langle V, V \rangle = 1\}$.
$S^2_1(1)$ is  a timelike surface in $\R^3_1$ known also  as the de Sitter space.

Assume that the curve $c: l = l(v) = l(w^1(v),w^2(v)), \, v \in J$  on
$S^2_1(1)$ is parameterized by the arc-length, i.e. $\langle l'(v), l'(v) \rangle = 1$.
Similarly to the elliptic case we consider an  orthonormal frame field
$\{ l(v), t(v), n(v)\}$ in $\R^3_1$, such that $t(v) = l'(v)$ and $\langle n(v), n(v) \rangle = -1$.
With respect to this  frame field we have the following decompositions of the vector fields $l'(v)$, $t'(v)$, $n'(v)$:
\begin{equation} \label{E:Eq8}
\begin{array}{l}
\vspace{2mm}
l' = t;\\
\vspace{2mm}
t' = - \kappa \,n - l;\\
\vspace{2mm} n' = - \kappa \,t,
\end{array}
\end{equation}
which can be considered as Frenet formulas of $c$ on $S^2_1(1)$.
The function $\kappa (v)= \langle t'(v), n(v) \rangle$ is the spherical
curvature of $c$ on $S^2_1(1)$.

We assume that $(f'(u))^2 + (g'(u))^2 = 1$.
Denote  $x = z_u = f' \,l + g'\,e_1, \,\, y = \ds{\frac{z_v}{f} = t}$ and
consider the orthonormal normal frame field defined by:
$$n_1 = g'(u)\,l(v) - f'(u) \, e_1; \qquad n_2 =  n(v).$$
Thus we obtain a frame field $\{x,y, n_1, n_2\}$ of $\mathcal{M}''_m$,
such that $\langle n_1, n_1 \rangle =1$, $\langle n_2, n_2 \rangle =
- 1$, $\langle n_1, n_2 \rangle =0$.

Using  \eqref{E:Eq8} we get
the following derivative formulas \cite{GM-new}:
\begin{equation} \label{E:Eq9}
\begin{array}{ll}
\vspace{2mm} \nabla'_xx = - \kappa_m\,n_1; & \qquad
\nabla'_x n_1 = \kappa_m \,x;\\
\vspace{2mm} \nabla'_xy = 0;  & \qquad
\nabla'_y n_1 = \ds{ \frac{g'}{f}\,y};\\
\vspace{2mm} \nabla'_yx = \ds{\frac{f'}{f}}\,y;  & \qquad
\nabla'_x n_2 = 0;\\
\vspace{2mm} \nabla'_yy = \ds{- \frac{f'}{f}\,x  -
\frac{g'}{f}\,n_1 - \frac{\kappa}{f} \, n_2}; & \qquad
\nabla'_y n_2 = \ds{  - \frac{\kappa}{f}\,y},
\end{array}
\end{equation}
where $\kappa_m (u)= f'(u) g''(u) - g'(u) f''(u)$ is the curvature of the meridian curve $m$, and
$\kappa = \kappa (v)$ is the spherical curvature of $c$.

\vskip 50mm
\section{Meridian surfaces of elliptic or hyperbolic type with harmonic Gauss map}

In the present section we give the classification of the meridian surfaces of elliptic or hyperbolic type with harmonic Gauss map.

Let $\mathcal{M}'_m$ and  $\mathcal{M}''_m$ be  meridian surfaces of elliptic and hyperbolic type, respectively, and $\{x,y,n_1,n_2\}$ be the frame field
of $\mathcal{M}'_m$ (resp. $\mathcal{M}''_m$) defined in Section \ref{S:Mer}. This frame field  generates the following frame of the
Grassmannian manifold: $$\{x \wedge y, x \wedge n_1, x \wedge n_2,
y \wedge n_1, y \wedge n_2, n_1 \wedge n_2\}.$$ The indefinite
inner product on the Grassmannian manifold is given by
\begin{equation}\notag
\langle e_{i_1} \wedge e_{i_2}, f_{j_1} \wedge f_{j_2}  \rangle =
\det \left( \langle e_{i_k}, f_{j_l}  \rangle  \right).
\end{equation}
Thus we have
\begin{equation} \notag
\begin{array}{lll}
\vspace{2mm}
\langle x \wedge y, x \wedge y  \rangle = 1; & \qquad \langle x \wedge n_1, x \wedge n_1  \rangle = 1; & \qquad \langle x \wedge n_2, x \wedge n_2  \rangle = - 1;\\
\langle y \wedge n_1, y \wedge n_1  \rangle = 1; & \qquad \langle y \wedge n_2, y \wedge n_2  \rangle = - 1; & \qquad \langle n_1 \wedge n_2, n_1 \wedge n_2  \rangle = - 1,
\end{array}
\end{equation}
and all other scalar products are equal to zero.

The Gauss map $G$ of $\mathcal{M}'_m$ (resp. $\mathcal{M}''_m$)  is defined by $G(p)=(x\wedge y)(p)$,
$p \in \mathcal{M}'_m$ (resp. $p \in \mathcal{M}''_m$). Then the Laplacian of the Gauss map is given by the
formula
\begin{equation}\label{Eq10}
\Delta G = - \nabla'_x\nabla'_x G + \nabla'_{\nabla_x x} G -
\nabla'_y\nabla'_y G + \nabla'_{\nabla_y y} G.
\end{equation}

Using \eqref{E:Eq5}, \eqref{E:Eq9} and \eqref{Eq10}, we obtain that in the elliptic case the Laplacian of the Gauss map
is expressed as
\begin{equation}\label{Eq11}
\Delta G = \frac{\kappa^2 - g'^ 2 - f^2 \kappa_m^2}{f^2}\, x\wedge y - \frac{\kappa'}{f^2}\, x\wedge n_1 -
\frac{\kappa f'}{f^2}\, y\wedge n_1 + \frac{f (f \kappa_m)'- f' g'}{f^2}\, y\wedge n_{2},
\end{equation}
 and in the hyperbolic case the Laplacian of the Gauss map
is given by
\begin{equation}\label{Eq12}
\Delta G = \frac{- \kappa^2 + g'^ 2 + f^2 \kappa_m^2}{f^2}\, x\wedge y + \frac{\kappa'}{f^2}\, x\wedge n_2
+ \frac{ f' g' - f (f \kappa_m)'}{f^2}\, y\wedge n_1 + \frac{\kappa f'}{f^2}\, y\wedge n_2,
\end{equation}
where $\kappa^{\prime }= \displaystyle{\frac{d}{dv}(\kappa)}$.

\begin{thm} \label{T:harmonic-ell}
Let $\mathcal{M}'_m$  be a meridian surface of elliptic type, defined by \eqref{E:Eq2}.
The Gauss map of $\mathcal{M}'_m$  is harmonic if and only if $\mathcal{M}'_m$  is part of a plane.
\end{thm}

\noindent {\it Proof:}
First, we suppose that the Gauss map of $\mathcal{M}'_m$ is harmonic, i.e. $\Delta G = 0$.
Then, from \eqref{Eq11} it follows that
\begin{equation} \notag
\begin{array}{lll}
\vspace{2mm}
\kappa^2 - g'^ 2 - f^2 \kappa_m^2 = 0;\\
\vspace{2mm}
\kappa' = 0;\\
\vspace{2mm}
\kappa f' = 0;\\
f (f \kappa_m)'- f' g' = 0.
\end{array}
\end{equation}
In the elliptic case we have $f'^2 \geq 1$, since $f'^2 - g'^2 = 1$. Hence, the above equalities imply
\begin{equation} \notag
\begin{array}{l}
\vspace{2mm}
\kappa = 0;\\
\vspace{2mm}
g' = 0;\\
\kappa_m = 0.
\end{array}
\end{equation}
Using \eqref{E:Eq5} we get that $\mathcal{M}'_m$ is totally geodesic, i.e. $\mathcal{M}'_m$ is part of a plane.

Conversely, if $\mathcal{M}'_m$  is totally geodesic, then $\Delta G = 0$.
\qed

\begin{thm} \label{T:harmonic-hyp}
Let $\mathcal{M}''_m$  be a meridian surface of hyperbolic type, defined by \eqref{E:Eq6}.
The Gauss map of $\mathcal{M}''_m$  is harmonic if and only if one of the following cases holds:

(i) $\mathcal{M}''_m$  is part of a plane;

(ii) the curve $c$ has spherical
curvature $\kappa = \pm 1$ and the meridian curve $m$ is determined by $f(u) = a; \, \,
g(u) = \pm u + b$, where $a=const$, $b=const$. In this case $\mathcal{M}''_m$  is a marginally trapped
developable ruled surface in $\E^4_1$.
\end{thm}

\noindent {\it Proof:}
Suppose that the Gauss map of $\mathcal{M}''_m$ is harmonic, i.e. $\Delta G = 0$.
Then, from \eqref{Eq12} it follows that
\begin{equation} \label{E:Eq13}
\begin{array}{lll}
\vspace{2mm}
\kappa^2 - g'^ 2 - f^2 \kappa_m^2 = 0;\\
\vspace{2mm}
\kappa' = 0;\\
\vspace{2mm}
\kappa f' = 0;\\
f (f \kappa_m)'- f' g' = 0.
\end{array}
\end{equation}
In the hyperbolic case we have $f'^2 \leq 1$, since $f'^2 + g'^2 = 1$.
Hence, from the third equality of \eqref{E:Eq13} we get the following two cases:

\vskip 1mm
Case (i):  $\kappa =0$. Then, the first equality of \eqref{E:Eq13} implies
$g' = 0;\,\, \kappa_m = 0$. Using \eqref{E:Eq9} we get that $\mathcal{M}''_m$ is totally geodesic, i.e. $\mathcal{M}''_m$ is part of a plane.

\vskip 1mm
Case (ii): $\kappa \neq 0$. Then $f' =0$, i.e. $f(u) = a = const$, and $g'^2 = 1$, i.e. $g(u) = \pm u + b$, $b = const$.
In this case $\kappa_m = 0$.
The second equality of \eqref{E:Eq13} implies $\kappa = const$. From the first equality of \eqref{E:Eq13} we obtain $\kappa^2 = g'^2$.
Hence, $\kappa = \varepsilon g'$, where $\varepsilon = \pm 1$. In this case  derivative formulas \eqref{E:Eq9} take the form:
\begin{equation} \label{E:Eq14}
\begin{array}{ll}
\vspace{2mm} \nabla'_xx = 0; & \qquad
\nabla'_x n_1 = 0;\\
\vspace{2mm} \nabla'_xy = 0;  & \qquad
\nabla'_y n_1 = \ds{\pm \frac{1}{a}\,y};\\
\vspace{2mm} \nabla'_yx = 0;  & \qquad
\nabla'_x n_2 = 0;\\
\vspace{2mm} \nabla'_yy = \ds{\mp  \frac{1}{a} ( n_1 + \varepsilon n_2)}; & \qquad
\nabla'_y n_2 = \ds{ \mp  \frac{\varepsilon}{a}\,y}.
\end{array}
\end{equation}
 $\mathcal{M}''_m$ is a ruled surface, since the meridian curve $m$ is a straight line.
 From \eqref{E:Eq14} we have that $\nabla'_x n_1 = 0; \,\, \nabla'_x n_2 = 0$, i.e. the normal space is constant at the points of each generator.
 Hence,  $\mathcal{M}''_m$ is developable.
 Moreover, the mean curvature vector field $H$ is given by
 $$H = \ds{\mp  \frac{1}{2a} ( n_1 + \varepsilon n_2)},$$
 which implies that $\langle H, H \rangle = 0$. Hence, $\mathcal{M}''_m$ is a marginally trapped surface.

\vskip 1mm
Conversely, if one of the cases (i) or (ii) holds, then  by straightforward calculations we get $\Delta G = 0$, i.e. $\mathcal{M}''_m$  has harmonic Gauss map.
\qed

\vskip 3mm
\noindent
\emph{Remark}: In the Euclidean space $\E^4$ planes are the only surfaces with harmonic Gauss map.
However, in the Minkowski space $\E^4_1$ there are surfaces with harmonic Gauss map which are not planes.
Theorem \ref{T:harmonic-ell} and Theorem  \ref{T:harmonic-hyp} show that in the class of the meridian surfaces of elliptic type
there are no surfaces with harmonic Gauss map other than  planes, while in the class of
the meridian surfaces of hyperbolic type we obtain surfaces with harmonic Gauss map, which are not planes.

\section{Meridian surfaces of elliptic or hyperbolic type with pointwise 1-type  Gauss map of first kind}

In this section we classify the meridian surfaces of elliptic or hyperbolic type with pointwise 1-type  Gauss map of first kind, i.e.
the Gauss map $G$ satisfies the condition
$$\Delta G = \lambda G$$
for some non-zero smooth function $\lambda$.

First, let us consider the meridian surface of elliptic type $\mathcal{M}'_m$, defined by \eqref{E:Eq2}.
So, the Laplacian of the Gauss map is given by formula  \eqref{Eq11}.

\begin{thm} \label{T:First kind - ell}
Let $\mathcal{M}'_m$ be a meridian surface of elliptic  type, defined by \eqref{E:Eq2}.
Then $\mathcal{M}'_m$ has pointwise 1-type Gauss map of first kind if and only if the curve $c$ has zero spherical
curvature and the meridian curve $m$ is determined by a solution $f(u)$ of the following differential
equation
\begin{equation} \label{E:Eq15}
f \left(\frac{f f''}{\sqrt{f'^2-1}}\right)' - f' \sqrt{f'^2-1} = 0,
\end{equation}
$g(u)$ is defined by $g'(u) = \sqrt{f'^2(u) - 1}$.
\end{thm}

\vskip 2mm
\noindent {\it Proof:}
From \eqref{Eq11} it follows that $\Delta G = \lambda G$ if and only if
\begin{equation} \label{E:Eq16}
\begin{array}{lll}
\vspace{2mm}
\kappa' = 0;\\
\vspace{2mm}
\kappa f' = 0;\\
f (f \kappa_m)'- f' g' = 0.
\end{array}
\end{equation}
Let $\mathcal{M}'_m$ be of pointwise 1-type Gauss map of first kind. Since $f' \neq 0$, from the second equality of  \eqref{E:Eq16} we get that $\kappa = 0$.
If we suppose that $g' = 0$, then $\kappa_m = 0$ and \eqref{Eq11} implies $\Delta G = 0$, which contradicts the assumption $\lambda \neq 0$. Hence, $g' \neq 0$. Then
using that $f'^2 - g'^2 = 1$ we obtain $\kappa_m = \ds{\frac{f''}{\sqrt{f'^2 - 1}}}$.
So, the third equality of \eqref{E:Eq16} takes the form \eqref{E:Eq15}.
Conversely, if $\kappa = 0$ and $f(u)$ is a solution of \eqref{E:Eq15}, then $\Delta G = \lambda G$.
\qed

\vskip 3mm
In the next theorem we give the classification of the meridian surfaces of hyperbolic type with pointwise 1-type Gauss map of first kind.

\begin{thm} \label{T:First kind - hyp}
Let $\mathcal{M}''_m$ be a meridian surface of hyperbolic type, defined by \eqref{E:Eq6}.
Then $\mathcal{M}''_m$ has pointwise 1-type Gauss map of first kind if and only if one of the following cases holds:

(i) the curve $c$ has zero spherical
curvature and the meridian curve $m$ is determined by a solution $f(u)$ of the following differential
equation
\begin{equation} \label{E:Eq17}
f \left(\frac{f f''}{\sqrt{1-f'^2}}\right)' + f' \sqrt{1-f'^2} = 0,
\end{equation}
$g(u)$ is defined by $g'(u) = \sqrt{1-f'^2(u)}$;

(ii) the curve $c$ has non-zero constant spherical
curvature $k$ ($\kappa \neq \pm 1$) and the meridian curve $m$ is determined by $f(u) = a; \, \,
g(u) = \pm u + b$, where $a=const$, $b=const$. Moreover, $\mathcal{M}''_m$  is a
developable ruled surface lying in a constant hyperplane $\E^3_1$ (if $\kappa^2-1>0$) or $\E^3$ (if $\kappa^2-1<0$) of $\E^4_1$.
\end{thm}

\vskip 2mm
\noindent {\it Proof:} Let  $\mathcal{M}''_m$  be a meridian surface of hyperbolic type, defined by \eqref{E:Eq6}.
So, the Laplacian of the Gauss map is given by formula  \eqref{Eq12}.
From \eqref{Eq12} it follows that $\Delta G = \lambda G$ if and only if
\begin{equation} \label{E:Eq18}
\begin{array}{lll}
\vspace{2mm}
\kappa' = 0;\\
\vspace{2mm}
\kappa f' = 0;\\
f' g' - f (f \kappa_m)' = 0.
\end{array}
\end{equation}
Let $\mathcal{M}''_m$ be of pointwise 1-type Gauss map of first kind.
From the second equality of \eqref{E:Eq18} we get the following two cases:

\vskip 1mm
Case (i):  $\kappa =0$. If we suppose that $g' = 0$, then $\kappa_m = 0$ and \eqref{Eq12} implies $\Delta G = 0$, which contradicts the assumption
 $\lambda \neq 0$. Hence, $g' \neq 0$. Then
using that $f'^2 + g'^2 = 1$ we obtain $\kappa_m = \ds{- \frac{f''}{\sqrt{1 - f'^2}}}$.
Thus, the third equality of \eqref{E:Eq18} takes the form \eqref{E:Eq17}.

\vskip 1mm
Case (ii):  $\kappa \neq 0$. Then $f' =0$, i.e. $f(u) = a = const$, and $g' = \pm 1$, i.e. $g(u) = \pm u + b$, $b = const$.
In this case $\kappa_m = 0$. The first equality of \eqref{E:Eq18} implies $\kappa = const$.
Then the Laplacian of the Gauss map takes the form:
\begin{equation} \notag
\Delta G = \frac{1 - \kappa^2}{a^2}\, x\wedge y,
\end{equation}
which implies that the surface  $\mathcal{M}''_m$   has 1-type Gauss map, since $\lambda =  \ds{\frac{1 - \kappa^2}{a^2}} = const$.
If  $\kappa^2 =1$ then $\Delta G = 0$, which contradicts the assumption
 $\lambda \neq 0$. Hence,  $\kappa^2 \neq 1$, i.e. $\kappa \neq \pm 1$.
In this case derivative formulas \eqref{E:Eq9} take the form:
\begin{equation} \label{Eq19}
\begin{array}{ll}
\vspace{2mm} \nabla'_xx = 0; & \qquad
\nabla'_x n_1 = 0;\\
\vspace{2mm} \nabla'_xy = 0;  & \qquad
\nabla'_y n_1 = \ds{\pm \frac{1}{a}\,y};\\
\vspace{2mm} \nabla'_yx = 0;  & \qquad
\nabla'_x n_2 = 0;\\
\vspace{2mm} \nabla'_yy = \ds{\mp  \frac{1}{a} n_1 - \frac{\kappa}{a} n_2}; & \qquad
\nabla'_y n_2 = \ds{ - \frac{\kappa}{a}\,y}.
\end{array}
\end{equation}
$\mathcal{M}''_m$ is a developable ruled surface, since the meridian curve $m$ is a straight line and
 $\nabla'_x n_1 = 0; \,\, \nabla'_x n_2 = 0$.
Now we shall prove that $\mathcal{M}''_m$  lies in a constant hyperplane of $\E^4_1$.

An arbitrary orthonormal frame  $\{n, n^{\bot}\}$ of the normal bundle is determined by
\begin{equation} \label{Eq20}
\begin{array}{l}
\vspace{2mm}
n = \cosh \theta \,n_1 + \sinh \theta \,n_2\\
\vspace{2mm}
n^{\bot} = \sinh \theta \,n_1 + \cosh \theta \,n_2,
\end{array}
\end{equation}
for some smooth function $\theta$. Note that $n$ is spacelike and $n^{\bot}$ is timelike.
Using \eqref{Eq19} and \eqref{Eq20} we get
\begin{equation} \label{Eq21}
\begin{array}{l}
\vspace{2mm}
\nabla'_x n = \theta'_u \, n^{\bot};\\
\vspace{2mm}
\nabla'_y n = \ds{\frac{\theta'_v}{a} \, n^{\bot} - \frac{1}{a} (\mp \cosh \theta + \kappa \sinh \theta)\, y};\\
\vspace{2mm}
\nabla'_x n^{\bot} = \theta'_u \, n;\\
\vspace{2mm}
\nabla'_y n^{\bot} = \ds{\frac{\theta'_v}{a} \, n - \frac{1}{a} (\mp \sinh \theta + \kappa \cosh \theta)\, y}.
\end{array}
\end{equation}

In the case $\kappa^2 -1 >0$ we choose $\theta = \ds{ \pm \frac{1}{2} \ln \left(\frac{\kappa + 1}{\kappa - 1}\right)}$. Then
$\mp \cosh \theta + \kappa \sinh \theta = 0$. Hence,  from \eqref{Eq19} and \eqref{Eq21}  it follows that
\begin{equation} \notag
\begin{array}{ll}
\vspace{2mm} \nabla'_xx = 0; & \qquad
\nabla'_x n = 0;\\
\vspace{2mm} \nabla'_xy = 0;  & \qquad
\nabla'_y n = 0;\\
\vspace{2mm} \nabla'_yx = 0;  & \qquad
\nabla'_x n^{\bot} = 0;\\
\vspace{2mm} \nabla'_yy = \ds{- \frac{\kappa^2 - 1}{a \kappa} \cosh \ln \left(\frac{\kappa + 1}{\kappa - 1}\right)^{\pm \frac{1}{2}} \,n^{\bot}}; & \qquad
\nabla'_y n^{\bot} = \ds{- \frac{\kappa^2 - 1}{a \kappa} \cosh \ln \left(\frac{\kappa + 1}{\kappa - 1}\right)^{\pm \frac{1}{2}} \,y}.
\end{array}
\end{equation}
The last equalities imply that $n = const$ and the surface  $\mathcal{M}''_m$  lies in the constant hyperplane $\E^3_1 = \span \{x, y, n^{\bot} \}$ of $\E^4_1$.

In the case $\kappa^2 -1 <0$ we choose $\theta = \ds{\pm \frac{1}{2} \ln \left(\frac{\kappa + 1}{1 - \kappa}\right)}$. Then
$\mp \sinh \theta + \kappa \cosh \theta = 0$. Hence, formulas \eqref{Eq19} and \eqref{Eq21}  imply that
\begin{equation} \notag
\begin{array}{ll}
\vspace{2mm} \nabla'_xx = 0; & \qquad
\nabla'_x n = 0;\\
\vspace{2mm} \nabla'_xy = 0;  & \qquad
\nabla'_y n = \ds{- \frac{\kappa^2 - 1}{a \kappa} \sinh \ln \left(\frac{\kappa + 1}{1 - \kappa}\right)^{\pm \frac{1}{2}} \,y};\\
\vspace{2mm} \nabla'_yx = 0;  & \qquad
\nabla'_x n^{\bot} = 0;\\
\vspace{2mm} \nabla'_yy = \ds{\frac{\kappa^2 - 1}{a \kappa} \sinh \ln \left(\frac{\kappa + 1}{1 - \kappa}\right)^{\pm \frac{1}{2}} \,n}; & \qquad
\nabla'_y n^{\bot} = 0.
\end{array}
\end{equation}
From the last equalities we get  that $n^{\bot} = const$ and the surface  $\mathcal{M}''_m$  lies in the constant hyperplane $\E^3 = \span \{x, y, n \}$ of $\E^4_1$.

\vskip 1mm
Conversely, if one of the cases (i) or (ii) holds, then by straightforward calculations   it can be seen that $\Delta G = \lambda G$, i.e. $\mathcal{M}''_m$
has pointwise 1-type Gauss map of first kind.
\qed

\section{Meridian surfaces of elliptic or hyperbolic type with pointwise 1-type  Gauss map of second kind}

In this section we give the classification of the meridian surfaces of elliptic or hyperbolic type with pointwise 1-type  Gauss map of second kind, i.e.
the Gauss map $G$ satisfies the condition
\begin{equation} \label{E:Eq22}
\Delta G = \lambda (G + C)
\end{equation}
for some non-zero smooth function $\lambda$ and a constant vector $C \neq 0$.

First we consider meridian surfaces of elliptic type with pointwise 1-type  Gauss map of second kind. They are classified by the following theorem.

\begin{thm} \label{T:Second kind - ell}
Let $\mathcal{M}'_m$ be a meridian surface of elliptic  type, defined by \eqref{E:Eq2}.
Then $\mathcal{M}'_m$ has pointwise 1-type Gauss map of second kind if and only if one of the following cases holds:

(i) the curve $c$ has non-zero constant spherical curvature $\kappa$ and the meridian curve $m$
is determined by $f(u) = \pm u + a; \, \, g(u) = b$, where $a=const$, $b=const $. In this case $\mathcal{M}'_m$  is  a developable ruled surface lying in
a constant hyperplane $\E^3$ of $\E^4_1$.

(ii) the curve $c$ has constant spherical curvature $\kappa$ and the meridian curve $m$
is determined by $f(u) = a u + a_1; \, \, g(u) = b u + b_1$, where $a$, $a_1$, $b$
and $b_1$ are constants, $a^2 \geq 1$, $a^2 - b^2 = 1$. In this case $\mathcal{M}'_m$   is either a
marginally trapped developable ruled surface (if $\kappa^2 = b^2$) or a developable ruled
surface lying in a constant hyperplane $\E^3$ (if $\kappa^2 - b^2 > 0$) or $\E^3_1$ (if $\kappa^2 - b^2 < 0$) of $\E^4_1$.

(iii)  the curve $c$ has zero  spherical curvature  and the meridian curve $m$ is determined by the solutions of the following differential
equation
\begin{equation*}
\left(\ln \frac{\sqrt{f'^2 - 1} \left( f (f'^2 - 1)(f f'')' - f^2 f' f''^2 - f' (f'^2 - 1)^2 \right)}{(f'^2 - 1)^2 + f^2 f''^2 - f f' (f'^2 - 1) (f f'')'}
\right)'= {\frac{f'f''}{f'^2 - 1}}.
\end{equation*}
$g(u)$ is defined by $g'(u) = \sqrt{f'^2(u) - 1}$.
\end{thm}

\vskip 2mm
\noindent {\it Proof:}
Let $\mathcal{M}'_m$ be a meridian surface of elliptic  type, defined by \eqref{E:Eq2}.
 Suppose that $\mathcal{M}'_m$ has pointwise 1-type Gauss map of second kind. Then equations \eqref{Eq11} and \eqref{E:Eq22} imply
\begin{equation}\label{E:Eq23}
\left( \frac{\kappa^2 - g'^ 2 - f^2 \kappa_m^2}{f^2} - \lambda \right)\, x\wedge y - \frac{\kappa'}{f^2}\, x\wedge n_1 -
\frac{\kappa f'}{f^2}\, y\wedge n_1 + \frac{f (f \kappa_m)'- f' g'}{f^2}\, y\wedge n_{2} = \lambda C.
\end{equation}
Since $\lambda \neq 0$, from \eqref{E:Eq23} we get
\begin{equation}\label{E:Eq24}
\begin{array}{l}
\vspace{2mm}
\langle C,  x\wedge y \rangle = \ds{\frac{\kappa^2 - g'^ 2 - f^2 \kappa_m^2}{\lambda f^2} - 1};\\
\vspace{2mm}
\langle C,  x\wedge n_1 \rangle = \ds{ - \frac{\kappa'}{\lambda f^2}};\\
\vspace{2mm}
\langle C,  y\wedge n_1 \rangle = \ds{- \frac{\kappa f'}{\lambda f^2}};\\
\vspace{2mm}
\langle C,  y\wedge n_2 \rangle = \ds{- \frac{f (f \kappa_m)'- f' g'}{\lambda f^2}};\\
\vspace{2mm}
\langle C,  x\wedge n_2 \rangle = 0;\\
\vspace{2mm}
\langle C,  n_1 \wedge n_2 \rangle = 0.
\end{array}
\end{equation}
Differentiating the last two equalities of \eqref{E:Eq24} with respect to $u$ and $v$ we obtain
\begin{equation}\label{E:Eq25-0}
\begin{array}{l}
\vspace{2mm}
\kappa_m \langle C,  x\wedge n_1 \rangle = 0;\\
\vspace{2mm}
g' \langle C,  x\wedge y \rangle + f' \langle C,  y\wedge n_2 \rangle  = 0;\\
\vspace{2mm}
g' \langle C,  y\wedge n_1 \rangle + \kappa \langle C,  y\wedge n_2 \rangle  = 0.
\end{array}
\end{equation}
Hence, equalities \eqref{E:Eq24} and \eqref{E:Eq25-0} imply
\begin{equation}\label{E:Eq25}
\begin{array}{l}
\vspace{2mm}
\kappa' \kappa_m =0;\\
\vspace{2mm}
\kappa (f \kappa_m)' = 0;\\
\vspace{2mm}
g'(1 + \kappa^2 - f^2 \kappa_m^2) - f f' (f \kappa_m)' = \lambda f^2 g'.
\end{array}
\end{equation}

We distinguish the following cases.

\vskip 2mm Case I: $g' = 0$. Then $\kappa \neq 0$ (otherwise the Gauss map is harmonic).
From \eqref{E:Eq23} we get
\begin{equation}\label{E:Eq26}
C = \left( \frac{\kappa^2 }{\lambda f^2} - 1\right)\, x\wedge y - \frac{\kappa'}{\lambda f^2}\, x\wedge n_1 -
\frac{\kappa f'}{\lambda f^2}\, y\wedge n_1.
\end{equation}
Using \eqref{E:Eq5} and \eqref{E:Eq26} we obtain
\begin{equation*}
\begin{array}{lll}
\vspace{2mm} \nabla'_xC & = & \displaystyle{\kappa^2 \left({\frac{1}{\lambda f^2}} \right)'_u \, x\wedge y - \kappa'\left( {\frac{1}{\lambda f^2}} \right)'_u\, x\wedge n_{1} -
\kappa \left( {\frac{f'}{\lambda f^2}} \right)'_u \,
y\wedge n_{1};} \\
\vspace{2mm} \nabla'_yC & = & \displaystyle{\frac{\kappa}{\lambda^2
f^3} \left( 3 \kappa'\lambda - \kappa \lambda'_v \right)
\, x\wedge y + \frac{1}{\lambda^2 f^3} \left(- \kappa''\lambda + k'\lambda'_v + \kappa^3 \lambda + \kappa
\lambda - \kappa \lambda^2 f^2 \right)\, x\wedge n_{1}} \\
&  & \displaystyle{+ \frac{f'}{\lambda^2 f^3} \left(- 2 \kappa'\lambda + \kappa \lambda'_v \right) \, y\wedge n_{1}}.
\end{array}
\end{equation*}
The last formulas imply that $C = const$ if and only if $\kappa = const$ and
$\lambda = \displaystyle{\frac{\kappa^2 + 1}{f^2}}$.
In this case
\begin{equation*}
\Delta G = {\frac{\kappa^2 + 1}{f^2}}(G + C),
\end{equation*}
where $C = \ds{-\frac{1}{\kappa^2 + 1} ( x\wedge y + \kappa f' \, y\wedge n_1)}$.
From $g' = 0$ it follows that $\kappa_m = 0$ and the meridian curve $m$ is determined by
$f(u) = \pm u + a; \, \, g(u) = b$, where $a=const$, $b=const $. The surface $\mathcal{M}'_m$  is a
developable ruled surface lying in the hyperplane $\E^3= \span\{x, y, n_1\}$, since  $\nabla'_x n_2 = 0; \, \, \nabla'_y n_2 =0$.

\vskip 2mm Case II: $g' \neq 0$. Then the last equality of \eqref{E:Eq25} implies
\begin{equation} \notag
\lambda = \frac{1}{g' f^2} \left( g' (1+\kappa^2 - f^2 \kappa_m^2) - f f' (f \kappa_m)'\right).
\end{equation}
It follows from the first two equalities of \eqref{E:Eq25} that
there are three subcases.

\vskip 2mm \hskip 10mm 1. $\kappa_m = 0$. In this subcase the Laplacian of $G$ is given by
\begin{equation}\label{E:Eq28}
\Delta G =  \frac{\kappa^2 - g'^ 2}{f^2} \, x\wedge y - \frac{\kappa'}{f^2}\, x\wedge n_1 -
\frac{\kappa f'}{f^2}\, y\wedge n_1 - \frac{f' g'}{f^2}\, y\wedge n_{2}
\end{equation}
and the function $\lambda$ is expressed as $\lambda = \ds{\frac{\kappa^2 + 1}{f^2}}$.
Now, equalities \eqref{E:Eq22} and \eqref{E:Eq28} imply
\begin{equation}\label{E:Eq29}
C= - \frac{1}{1+\kappa^2} \left(f'^2 \, x\wedge y + \kappa'\, x\wedge n_1 + \kappa f' \, y\wedge n_1 + f' g'\, y\wedge n_{2}\right).
\end{equation}
Using  formulas \eqref{E:Eq5} in the case $\kappa_m = 0$ and \eqref{E:Eq29} we obtain
\begin{equation*}
\begin{array}{lll}
\vspace{2mm} \nabla'_xC & = & \ds{- \frac{1}{1+\kappa^2}  \left( 2f' f'' \, x\wedge y + \kappa f'' \, y\wedge n_1 + (f' g'' +g' f'')\, y\wedge n_{2}\right)}; \\
\vspace{2mm} \nabla'_yC & = & \ds{\frac{\kappa \kappa'}{f(1+\kappa^2)^2} \left( 2f'^2 + 1+\kappa^2\right) \, x\wedge y
+ \frac{1}{f(1+\kappa^2)^2} \left(2 \kappa \kappa'^2 - \kappa''(1 + \kappa^2) \right)\, x\wedge n_{1}} \\
&  & \ds{- \frac{2 \kappa' f'}{f(1+\kappa^2)^2} \, y\wedge n_1} + \ds{\frac{2 \kappa \kappa' f' g'}{f(1 + \kappa^2)^2} \, y\wedge n_2}.
\end{array}
\end{equation*}
The last formulas imply that $C = const$ if and only if $\kappa = const$ and $f'' =0$. Hence,  the meridian curve $m$ is determined by
$f(u) = a u + a_1; \, \, g(u) = b u + b_1$, where $a$, $a_1$, $b$ and $b_1$ are constants, $a^2 \geq 1$, $a^2 - b^2 = 1$.
Hence,  $\mathcal{M}'_m$  is a developable ruled surface, since  $\nabla'_x n_1 = 0; \, \, \nabla'_x n_2 =0$.
We shall prove that in the case $\kappa^2 = b^2$ the surface $\mathcal{M}'_m$ is marginally trapped and in the case $\kappa^2 \neq b^2$ the surface $\mathcal{M}'_m$ lies  in a constant hyperplane $\E^3$ or $\E^3_1$  of $\E^4_1$.
Indeed, if $\kappa = \varepsilon b$, $\varepsilon = \pm 1$, then from \eqref{E:Eq5} we get that $H = \ds{\frac{b}{2f}(\varepsilon n_1 + n_2)}$ and hence, $\langle H,H \rangle = 0$,
which implies that $\mathcal{M}'_m$ is a marginally trapped surface.
In the case $\kappa^2 - b^2 \neq 0$ we consider an orthonormal frame  $\{n, n^{\bot}\}$ of the normal bundle which is determined by equalities
\eqref{Eq20} for some function $\theta$. Hence, the derivatives of $n$ and $n^{\bot}$ satisfy
\begin{equation} \notag
\begin{array}{l}
\vspace{2mm}
\nabla'_x n = \theta'_u \, n^{\bot};\\
\vspace{2mm}
\nabla'_y n = \ds{\frac{\theta'_v}{f} \, n^{\bot} + \frac{1}{f} (b \sinh \theta - \kappa \cosh \theta)\, y};\\
\vspace{2mm}
\nabla'_x n^{\bot} = \theta'_u \, n;\\
\vspace{2mm}
\nabla'_y n^{\bot} = \ds{\frac{\theta'_v}{f} \, n + \frac{1}{f} (b\cosh \theta - \kappa \sinh \theta)\, y}.
\end{array}
\end{equation}
In the case $\kappa^2 - b^2 >0$ we choose $\theta = \ds{\frac{1}{2}\ln \left(\frac{\kappa + b}{\kappa - b}\right)}$. Then
$b\cosh \theta - \kappa \sinh \theta = 0$ and $\nabla'_x n^{\bot} = 0$, $\nabla'_y n^{\bot} = 0$.
In this case the surface  $\mathcal{M}'_m$  lies in the constant hyperplane $\E^3 = \span \{x, y, n \}$ of $\E^4_1$.
In the case $\kappa^2 - b^2 <0$ we choose $\theta = \ds{\frac{1}{2} \ln \left(\frac{b + \kappa}{b - \kappa}\right)}$. Then
$b \sinh \theta - \kappa \cosh \theta = 0$ and $\nabla'_x n = 0$, $\nabla'_y n = 0$.
In this case the surface  $\mathcal{M}'_m$  lies in the constant hyperplane $\E^3_1 = \span \{x, y, n^{\bot} \}$ of $\E^4_1$.

\vskip 2mm \hskip 10mm 2. $\kappa = 0$. In this subcase the Laplacian of $G$ is given by
\begin{equation}\label{E:Eq30}
\Delta G =  - \frac{g'^ 2 + f^2 \kappa_m^2}{f^2} \, x\wedge y + \frac{f(f\kappa_m)' - f' g'}{f^2}\, y\wedge n_{2}
\end{equation}
and the function $\lambda$ is expressed as $\lambda = \ds{\frac{1}{g' f^2} \left( g' (1 - f^2 \kappa_m^2) - f f' (f \kappa_m)'\right)}$.
Hence,  from equalities \eqref{E:Eq22} and \eqref{E:Eq30} we get
\begin{equation}\label{E:Eq31}
C = - \left( \frac{g'^ 2 + f^2 \kappa_m^2}{\lambda f^2} +1 \right)\, x\wedge y + \frac{f (f \kappa_m)'- f' g'}{\lambda f^2}\, y\wedge n_{2}.
\end{equation}
Denote $\psi = \ds{- \frac{g'^ 2 + f^2 \kappa_m^2}{\lambda f^2} - 1}$;   $\varphi = \ds{ \frac{f (f \kappa_m)'- f' g'}{\lambda f^2}}$. Then
$C = \psi \, x\wedge y + \varphi \, y\wedge n_{2}$.
Using  formulas \eqref{E:Eq5} in the case $\kappa = 0$ and \eqref{E:Eq31} we obtain
\begin{equation} \label{E:Eq32}
\begin{array}{l}
\vspace{2mm}
\nabla'_xC = (\psi' - \varphi \kappa_m)\, x\wedge y + (\varphi' - \psi \kappa_m)\, y\wedge n_{2}; \\
\vspace{2mm}
\nabla'_yC = 0.
\end{array}
\end{equation}
Using the expression of $\lambda$ we calculate that $\psi = \ds{\varphi \, \frac{f'}{g'}}$;
$\psi' - \varphi \kappa_m = \ds{\frac{f'}{g'}(\varphi' - \psi \kappa_m)}$.
Hence, formulas \eqref{E:Eq32} take the form
\begin{equation*}
\begin{array}{l}
\vspace{2mm}
\nabla'_xC = \ds{\frac{f'}{g'}(\varphi' - \psi \kappa_m)\, x\wedge y + (\varphi' - \psi \kappa_m)\, y\wedge n_{2}}; \\
\vspace{2mm}
\nabla'_yC = 0.
\end{array}
\end{equation*}
The last formulas imply that $C = const$ if and only if
\begin{equation} \label{E:Eq33}
\left(\ln \varphi \right)' = \ds{\frac{f'}{g'}\kappa_m}.
\end{equation}
Using that $f \kappa_m = \ds{\frac{f f''}{\sqrt{f'^2 - 1}}}$, we get
\begin{equation} \label{E:Eq34}
\varphi = \ds{\frac{\sqrt{f'^2 - 1} \left( f (f'^2 - 1) (f f'')' - f^2 f' f''^2 - f' (f'^2 - 1)^2 \right)}{(f'^2 - 1)^2 + f^2 f''^2 - f f' (f'^2 - 1) (f f'')'}}.
\end{equation}
Now, formulas \eqref{E:Eq33} and \eqref{E:Eq34} imply that $C=const$ if and
only if the function $f(u)$ is a solution of the following differential
equation
\begin{equation} \notag
\left( \ln \ds{\frac{\sqrt{f'^2 - 1} \left( f (f'^2 - 1) (f f'')' - f^2 f' f''^2 - f' (f'^2 - 1)^2 \right)}{(f'^2 - 1)^2 + f^2 f''^2 - f f' (f'^2 - 1) (f f'')'}} \right)' =
\ds{\frac{f' f''}{f'^2 - 1}}.
\end{equation}

\vskip2mm \hskip10mm 3. $\kappa = const \neq 0$ and $f \kappa_m = a = const$, $ a \neq 0$.
In this subcase the Laplacian of $G$ is given by
\begin{equation}\label{E:Eq36}
\Delta G =  \frac{\kappa^2 - a^2 - g'^ 2}{f^2} \, x\wedge y - \frac{\kappa f'}{f^2}\, y\wedge n_1 - \frac{f' g'}{f^2}\, y\wedge n_2
\end{equation}
and the function $\lambda$ is expressed as $\lambda = \ds{\frac{1 + \kappa^2 - a^2}{f^2}}$. Since $\lambda \neq 0$, we get $a^2 \neq 1 + \kappa^2$.
Now, equalities \eqref{E:Eq22} and \eqref{E:Eq36} imply
\begin{equation*}
C = \frac{1}{1 + \kappa^2 - a^2} \left( - f'^2\, x\wedge y  - \kappa f'\, y\wedge n_1 - f' g'\, y\wedge n_2 \right).
\end{equation*}
Then the derivatives of $C$ are expressed as
\begin{equation}\label{E:Eq37}
\begin{array}{l}
\vspace{2mm}
\nabla'_xC = \ds{- \frac{1}{1 + \kappa^2 - a^2} \left(f' f''\, x\wedge y  + \kappa f''\, y\wedge n_1 + g' f'' \, y\wedge n_2 \right)};\\
\vspace{2mm}
\nabla'_yC = 0.
\end{array}\end{equation}
It follows from formulas \eqref{E:Eq37}  that $C = const$ if and only if $f'' = 0$. But the condition  $f'' = 0$ implies $\kappa_m = 0$, which
contradicts the assumption that $f \kappa_m \neq 0$.

Consequently, in this subcase  there are no meridian surfaces of elliptic type with pointwise 1-type Gauss
map of second kind.

\vskip 2mm
Conversely, if one of the cases (i), (ii) or (iii) holds, then  it can easily be seen that $\Delta G = \lambda (G + C)$, i.e. $\mathcal{M}'_m$
has pointwise 1-type Gauss map of second kind.
\qed

\vskip 3mm
The next theorem gives the classification of the meridian surfaces of hyperbolic type with pointwise 1-type  Gauss map of second kind.

\begin{thm} \label{T:Second kind - hyp}
Let $\mathcal{M}''_m$ be a meridian surface of hyperbolic  type, defined by \eqref{E:Eq6}.
Then $\mathcal{M}''_m$ has pointwise 1-type Gauss map of second kind if and only if one of the following cases holds:

(i) the curve $c$ has non-zero constant spherical curvature $\kappa \neq \pm 1$ and the meridian curve $m$
is determined by $f(u) = \pm u + a; \, \, g(u) = b$, where $a=const$, $b=const $. In this case $\mathcal{M}''_m$  is  a developable ruled surface lying in
a constant hyperplane $\E^3_1$ of $\E^4_1$.

(ii) the curve $c$ has constant spherical curvature $\kappa$ and the meridian curve $m$
is determined by $f(u) = a u + a_1; \, \, g(u) = b u + b_1$, where $a$, $a_1$, $b$
and $b_1$ are constants,  $a^2 + b^2 = 1$. In this case $\mathcal{M}''_m$   is either a
marginally trapped developable ruled surface (if $\kappa^2 = b^2$) or a developable ruled surface lying in a constant hyperplane
$\E^3_1$ (if $\kappa^2 - b^2 > 0$) or $\E^3$ (if $\kappa^2 - b^2 < 0$) of $\E^4_1$.

(iii)  the curve $c$ has zero  spherical curvature  and the meridian curve $m$ is determined by the solutions of the following differential
equation
\begin{equation*}
\left(\ln \frac{\sqrt{1 - f'^2} \left( f (1 - f'^2)(f f'')' + f^2 f' f''^2 + f' (1 - f'^2)^2 \right)}{(1 - f'^2)^2 + f^2 f''^2 + f f' (1 - f'^2) (f f'')'}
\right)'= - \frac{f'f''}{1 - f'^2}.
\end{equation*}
$g(u)$ is defined by $g'(u) = \sqrt{1 - f'^2(u)}$.
\end{thm}

\vskip 2mm
\noindent {\it Proof:}
Let $\mathcal{M}''_m$ be a meridian surface of hyperbolic  type, defined by \eqref{E:Eq6}.
 Suppose that $\mathcal{M}''_m$ has pointwise 1-type Gauss map of second kind.
Then equations \eqref{Eq12} and \eqref{E:Eq22} imply
\begin{equation}\label{E:Eq23-h}
\left( \frac{ g'^ 2 - \kappa^2 + f^2 \kappa_m^2}{f^2} - \lambda \right)\, x\wedge y + \frac{\kappa'}{f^2}\, x\wedge n_2
+ \frac{f' g' - f (f \kappa_m)'}{f^2}\, y\wedge n_1 + \frac{\kappa f'}{f^2}\, y\wedge n_2  = \lambda C.
\end{equation}
Since $\lambda \neq 0$, from \eqref{E:Eq23-h} we get
\begin{equation}\label{E:Eq24-h}
\begin{array}{l}
\vspace{2mm}
\langle C,  x\wedge y \rangle = \ds{\frac{g'^ 2 - \kappa^2 + f^2 \kappa_m^2}{\lambda f^2} - 1};\\
\vspace{2mm}
\langle C,  x\wedge n_2 \rangle = \ds{ - \frac{\kappa'}{\lambda f^2}};\\
\vspace{2mm}
\langle C,  y\wedge n_1 \rangle = \ds{\frac{f' g' - f (f \kappa_m)'}{\lambda f^2}};\\
\vspace{2mm}
\langle C,  y\wedge n_2 \rangle = \ds{- \frac{\kappa f'}{\lambda f^2}};\\
\vspace{2mm}
\langle C,  x\wedge n_1 \rangle = 0;\\
\vspace{2mm}
\langle C,  n_1 \wedge n_2 \rangle = 0.
\end{array}
\end{equation}
Differentiating the last two equalities of \eqref{E:Eq24-h} with respect to $u$ and $v$  we obtain
\begin{equation}\label{E:Eq25-h}
\begin{array}{l}
\vspace{2mm}
\kappa' \kappa_m =0;\\
\vspace{2mm}
\kappa (f \kappa_m)' = 0;\\
\vspace{2mm}
g'(1 - \kappa^2 + f^2 \kappa_m^2) - f f' (f \kappa_m)' = \lambda f^2 g'.
\end{array}
\end{equation}

Similarly to the elliptic case, we distinguish the following cases.

\vskip 2mm Case I: $g' = 0$. Then $\kappa \neq 0$ (otherwise the Gauss map is harmonic).
From \eqref{E:Eq23-h} we get
\begin{equation}\notag
C = - \left( \frac{\kappa^2 }{\lambda f^2} + 1\right)\, x\wedge y + \frac{\kappa'}{\lambda f^2}\, x\wedge n_2 +
\frac{\kappa f'}{\lambda f^2}\, y\wedge n_2,
\end{equation}
which implies
\begin{equation*}
\begin{array}{lll}
\vspace{2mm}
\nabla'_xC & = & \displaystyle{- \kappa^2 \left({\frac{1}{\lambda f^2}} \right)'_u \, x\wedge y +
\kappa'\left( {\frac{1}{\lambda f^2}} \right)'_u\, x\wedge n_2 + \kappa \left( {\frac{f'}{\lambda f^2}} \right)'_u \, y\wedge n_2;} \\
\vspace{2mm}
\nabla'_yC & = & \displaystyle{- \frac{\kappa}{\lambda^2 f^3} \left( 3 \kappa'\lambda - \kappa \lambda'_v \right)
\, x\wedge y + \frac{1}{\lambda^2 f^3} \left(\kappa''\lambda - k'\lambda'_v + \kappa^3 \lambda - \kappa
\lambda + \kappa \lambda^2 f^2 \right)\, x\wedge n_2} \\
&  & \displaystyle{+ \frac{f'}{\lambda^2 f^3} \left(2 \kappa'\lambda - \kappa \lambda'_v \right) \, y\wedge n_2}.
\end{array}
\end{equation*}
It follows from the last formulas that $C = const$ if and only if $\kappa = const$ and
$\lambda = \displaystyle{\frac{1 - \kappa^2}{f^2}}$.
Since $\lambda \neq 0$ we get $\kappa \neq \pm 1$. The Laplacian of the Gauss map is expressed as
\begin{equation*}
\Delta G = {\frac{1 - \kappa^2}{f^2}}(G + C),
\end{equation*}
where $C = \ds{\frac{1}{\kappa^2 - 1} ( x\wedge y - \kappa f' \, y\wedge n_2)}$.
The condition  $g' = 0$ implies that $\kappa_m = 0$ and the meridian curve $m$ is determined by
$f(u) = \pm u + a; \, \, g(u) = b$, where $a=const$, $b=const $. The surface $\mathcal{M}''_m$  is a
developable ruled surface lying in the hyperplane $\E^3_1= \span\{x, y, n_2\}$.

\vskip 2mm Case II: $g' \neq 0$. Then the last equality of \eqref{E:Eq25-h} implies
\begin{equation} \notag
\lambda = \frac{1}{g' f^2} \left( g' (1-\kappa^2 + f^2 \kappa_m^2) - f f' (f \kappa_m)'\right).
\end{equation}
Similarly to the elliptic case we have to consider the following three subcases.

\vskip 2mm \hskip 10mm 1. $\kappa_m = 0$. In this subcase $\lambda = \ds{\frac{1 - \kappa^2}{f^2}}$, $\kappa \neq \pm 1$ and the Laplacian of $G$ is given by
\begin{equation}\label{E:Eq28-h}
\Delta G =  \frac{g'^ 2 - \kappa^2}{f^2} \, x\wedge y + \frac{\kappa'}{f^2}\, x\wedge n_2 + \frac{f' g'}{f^2}\, y\wedge n_1 +
\frac{\kappa f'}{f^2}\, y\wedge n_2.
\end{equation}
So, equalities \eqref{E:Eq22} and \eqref{E:Eq28-h} imply
\begin{equation}\label{E:Eq29-h}
C= \frac{1}{1 - \kappa^2} \left(- f'^2 \, x\wedge y + \kappa'\, x\wedge n_2  + f' g'\, y\wedge n_1 + \kappa f' \, y\wedge n_2\right).
\end{equation}
Using  \eqref{E:Eq5} and \eqref{E:Eq29-h} we obtain
\begin{equation*}
\begin{array}{lll}
\vspace{2mm} \nabla'_xC & = & \ds{ \frac{1}{1 - \kappa^2}  \left( - 2f' f'' \, x\wedge y + (f' g'' +g' f'')\, y\wedge n_1 + \kappa f'' \, y\wedge n_2 \right)}; \\
\vspace{2mm} \nabla'_yC & = & \ds{\frac{- \kappa \kappa'}{f(1-\kappa^2)^2} \left( 2f'^2 + 1 - \kappa^2\right) \, x\wedge y
+ \frac{1}{f(1 - \kappa^2)^2} \left(2 \kappa \kappa'^2 + \kappa''(1 - \kappa^2) \right)\, x\wedge n_2} \\
&  & \ds{+  \ds{\frac{2 \kappa \kappa' f' g'}{f(1 - \kappa^2)^2} \, y\wedge n_1} + \frac{2 \kappa' f'}{f(1 - \kappa^2)^2} \, y\wedge n_2}.
\end{array}
\end{equation*}
From the last formulas we get that $C = const$ if and only if $\kappa = const$ and $f'' =0$. Then the meridian curve $m$ is determined by
$f(u) = a u + a_1; \, \, g(u) = b u + b_1$, where $a$, $a_1$, $b$ and $b_1$ are constants, $a^2 + b^2 = 1$.
Hence,  $\mathcal{M}''_m$  is a developable ruled surface, since  $\nabla'_x n_1 = 0; \, \, \nabla'_x n_2 =0$.
Analogously to the elliptic case we prove that if $\kappa^2 = b^2$ then $\mathcal{M}''_m$ is  a marginally trapped surface, and if $\kappa^2-b^2 \neq 0$, then $\mathcal{M}''_m$
  lies  in a constant hyperplane $\E^3$ or $\E^3_1$  of $\E^4_1$.
Indeed, in the case $\kappa = \varepsilon b$, $\varepsilon =\pm 1$, from \eqref{E:Eq9} we obtain
$H= \ds{- \frac{b}{2f} (n_1 + \varepsilon n_2)}$, which implies that $\langle H, H\rangle = 0$, i.e. $\mathcal{M}''_m$ is  a marginally trapped surface.
If $\kappa^2 - b^2 >0$ we choose $\theta = \ds{\frac{1}{2}\ln \left(\frac{\kappa + b}{\kappa - b}\right)}$ and  find a suitable normal frame field  $\{n, n^{\bot}\}$ such that
 $\mathcal{M}''_m$  lies in the constant hyperplane $\E^3_1 = \span \{x, y, n^{\bot} \}$ of $\E^4_1$.
If $\kappa^2 - b^2 <0$ we choose $\theta = \ds{\frac{1}{2} \ln \left(\frac{b + \kappa}{b - \kappa}\right)}$ and get that
$\mathcal{M}''_m$  lies in the constant hyperplane $\E^3 = \span \{x, y, n \}$.

\vskip 2mm \hskip 10mm 2. $\kappa = 0$. In this subcase the Laplacian of $G$ is given by
\begin{equation}\label{E:Eq30-h}
\Delta G =  \frac{g'^ 2 + f^2 \kappa_m^2}{f^2} \, x\wedge y + \frac{f' g' - f(f\kappa_m)'}{f^2}\, y\wedge n_1
\end{equation}
and the function $\lambda$ is expressed as $\lambda = \ds{\frac{1}{g' f^2} \left( g' (1 + f^2 \kappa_m^2) - f f' (f \kappa_m)'\right)}$.
Hence,  equalities \eqref{E:Eq22} and \eqref{E:Eq30-h} imply that
\begin{equation}\notag
C = \left( \frac{g'^ 2 + f^2 \kappa_m^2}{\lambda f^2} - 1 \right)\, x\wedge y + \frac{f' g' - f (f \kappa_m)'}{\lambda f^2}\, y\wedge n_1.
\end{equation}
Denoting $\psi = \ds{\frac{g'^ 2 + f^2 \kappa_m^2}{\lambda f^2} - 1}$;   $\varphi = \ds{ \frac{f' g' - f (f \kappa_m)'}{\lambda f^2}}$,
as in the elliptic case we obtain that
$C = const$ if and only if
\begin{equation} \notag
\left(\ln \varphi \right)' = \ds{\frac{f'}{g'}\kappa_m}.
\end{equation}
In the hyperbolic case we have  $f \kappa_m = \ds{- \frac{f f''}{\sqrt{1 - f'^2}}}$ and the function $\varphi$ is expressed as follows:
\begin{equation} \notag
\varphi = \ds{\frac{\sqrt{1 - f'^2} \left( f (1 - f'^2) (f f'')' + f^2 f' f''^2 + f' (1 - f'^2)^2 \right)}{(1 - f'^2)^2 + f^2 f''^2 + f f' (1 - f'^2) (f f'')'}}.
\end{equation}
Consequently,  $C=const$ if and
only if the function $f(u)$ is a solution of the following differential
equation
\begin{equation} \notag
\left( \ln \ds{\frac{\sqrt{1 - f'^2} \left( f (1 - f'^2) (f f'')' + f^2 f' f''^2 + f' (1 - f'^2)^2 \right)}{(1 - f'^2)^2 + f^2 f''^2 + f f' (1 - f'^2) (f f'')'}} \right)' =
\ds{- \frac{f' f''}{1 - f'^2}}.
\end{equation}

\vskip2mm \hskip10mm 3. $\kappa = const \neq 0$ and $f \kappa_m = a = const$, $ a \neq 0$.
In this subcase we have
\begin{equation}\label{E:Eq36-h}
\Delta G =  \frac{g'^ 2 - \kappa^2 + a^2 }{f^2} \, x\wedge y + \frac{f' g'}{f^2}\, y\wedge n_1 + \frac{\kappa f'}{f^2}\, y\wedge n_2
\end{equation}
and $\lambda = \ds{\frac{1 - \kappa^2 + a^2}{f^2}}$. Since $\lambda \neq 0$ we get $a^2 \neq \kappa^2 - 1$.
Equalities \eqref{E:Eq22} and \eqref{E:Eq36-h} imply
\begin{equation*}
C = \frac{1}{1 - \kappa^2 + a^2} \left( - f'^2\, x\wedge y  + f' g'\, y\wedge n_1  + \kappa f'\, y\wedge n_2 \right).
\end{equation*}
Then the derivatives of $C$ are expressed as
\begin{equation}\label{E:Eq37-h}
\begin{array}{l}
\vspace{2mm}
\nabla'_xC = \ds{\frac{1}{1 -\kappa^2 + a^2} \left( - f' f''\, x\wedge y  + g' f'' \, y\wedge n_1 + \kappa f''\, y\wedge n_2  \right)};\\
\vspace{2mm}
\nabla'_yC = 0.
\end{array}\end{equation}
Formulas \eqref{E:Eq37-h} imply that $C = const$ if and only if $f'' = 0$. But the condition $f'' = 0$ implies $\kappa_m = 0$, which
contradicts the assumption  $f \kappa_m \neq 0$.

Consequently, in this subcase  there are no meridian surfaces of hyperbolic type with pointwise 1-type Gauss
map of second kind.

\vskip 2mm
Conversely, if one of the cases (i), (ii) or (iii) holds, then  it can be seen that $\Delta G = \lambda (G + C)$, i.e. $\mathcal{M}''_m$
has pointwise 1-type Gauss map of second kind.
\qed

\vskip 3mm
Meridian surfaces of parabolic type in $\E^4_1$ are defined as one-parameter systems of meridians of the rotational hypersurface with lightlike axis
analogously to the meridian surfaces of elliptic  and hyperbolic type \cite{GM7}.
Similarly to the elliptic and hyperbolic type one can classify the  meridian surfaces of parabolic type with pointwise 1-type Gauss map.

\vskip 3mm \textbf{Acknowledgements:} This paper is prepared
during the first named author's visit to the Institute of
Mathematics and Informatics at the Bulgarian Academy of Sciences,
Sofia, Bulgaria  in September 2014.


\begin{thebibliography}{99}

\bibitem{ABBO}
K. Arslan,  B. Bayram, B. Bulca, G. Ozturk, Meridian surfaces of Weingarten type in 4-dimensional Euclidean space $\E^4$, preprint available at
ArXiv:1305.3155.

\bibitem{ABM}
K. Arslan, B. Bulca, V. Milousheva, Meridian surfaces in
$\E^4$ with pointwise 1-type Gauss map,\emph{ Bull. Korean Math. Soc.} 51
(2014), no. 3, 911--922.

\bibitem{BB} C. Baikoussis, D. E. Blair, On the Gauss map of ruled
surfaces, \emph{Glasgow Math. J.} 34 (1992), no. 3, 355--359.

\bibitem{BCV} C. Baikoussis, B.-Y. Chen,  L. Verstraelen, Ruled
surfaces and tubes with finite type Gauss map, \emph{Tokyo J. Math.} 16
(1993), no. 2, 341--349.

\bibitem{BV} C. Baikoussis, L. Verstraelen, On the Gauss map of
helicoidal surfaces, \emph{Rend. Sem. Mat. Messina} Ser. II 2 (16)
(1993), 31--42.

\bibitem{Ch1} B.-Y. Chen, \emph{Total mean curvature and submanifolds of
finite type}, Series in Pure Mathematics, 1. World Scientific
Publishing Co., Singapore, 1984.

\bibitem{Chen2}
Chen B.-Y., \emph{Pseudo-Riemannian geometry, $\delta$-invariants and applications},
World Scientific Publishing Co. Pte. Ltd., Hackensack, NJ, 2011.

\bibitem{Ch2} B.-Y. Chen, \emph{Finite type submanifolds and generalizations}, Universit\'{a} degli Studi di Roma ''La Sapienza'', Dipartimento
di Matematica IV, Rome, 1985.

\bibitem{Ch3} B.-Y. Chen, A report on submanifolds of finite type,
\emph{Soochow J. Math.} 22 (1996), no. 2, 117--337.

\bibitem{CCK} B.-Y. Chen, M. Choi, Y. H. Kim, Surfaces of
revolution with pointwise 1-type Gauss map, \emph{J. Korean Math. Soc.}
42 (2005), no. 3, 447--455.

\bibitem{CP} B.-Y. Chen, P. Piccinni, Submanifolds with finite type
Gauss map, \emph{Bull. Austral. Math. Soc.} 35 (1987), no. 2, 161--186.

\bibitem{CKY}
M. Choi, Y. H. Kim, Young, D. W. Yoon, Classification of
ruled surfaces with pointwise 1-type Gauss map in Minkowski
3-space, \emph{Taiwanese J. Math.} 15 (2011), no. 3, 1141--1161.

\bibitem{GM2}
G. Ganchev, Milousheva V., Invariants and Bonnet-type theorem
for  surfaces in $\R^4$, \emph{Cent. Eur. J. Math.} 8, no. 6  (2010), 993--1008.

\bibitem{GM6}
G. Ganchev,  Milousheva V.,  An invariant theory of marginally
trapped surfaces in the four-dimensional Minkowski space, \emph{J. Math. Phys.} 53  (2012), Article ID: 033705, 15 pp.

\bibitem{GM7}
G. Ganchev,  Milousheva V., Marginally trapped meridian surfaces of parabolic type in the four-dimensional Minkowski space,
\emph{Int. J. Geom. Methods  Mod. Phys.}, 10, no. 10  (2013), Article ID: 1350060, 17 pp.

\bibitem{GM-new}
G. Ganchev, Milousheva V., Meridian surfaces of elliptic or
hyperbolic type in the four-dimensional Minkowski space, preprint
available at ArXiv:1402.6112.

\bibitem{KY1} Y. H. Kim, D. W. Yoon, Ruled surfaces with finite type Gauss map in Minkowski spaces, \emph{Soochow J. Math.} 26 (2000), no. 1, 85--96.

\bibitem{KY2} Y. H. Kim, D. W. Yoon, Ruled surfaces with pointwise
1-type Gauss map, \emph{J. Geom. Phys.} 34 (2000), no. 3-4, 191--205.

\bibitem{KY2-a}
 Y. H. Kim, D. W. Yoon, Classifications of rotation surfaces in pseudo-Euclidean space.
\emph{ J. Korean Math. Soc.} 41 (2004), no. 2, 379–-396.

\bibitem{KY3} Y. H. Kim, D. W. Yoon, On the Gauss map of ruled
surfaces in Minkowski space, \emph{Rocky Mountain J. Math.} 35 (2005),
no. 5, 1555--1581.

\end{thebibliography}
\end{document}